\documentclass{amsart}

\UseRawInputEncoding

\usepackage{amsmath, amsthm}
\usepackage{amssymb, latexsym}
\usepackage{amsfonts}

 \usepackage{graphicx}
\usepackage[dvipsnames]{xcolor}

\usepackage{isomath} 
\usepackage{bm}

\newtheorem{theorem}[equation]{Theorem}

\newtheorem{lemma}[equation]{Lemma}

\newtheorem{definition}[equation]{Definition}

\newtheorem{corollary}[equation]{Corollary}

\newtheorem{Remark}[equation]{Remark}

\newtheorem{coro}[equation]{Corollary}
\newtheorem{prop}[equation]{Proposition}

\newtheorem{defi}[equation]{Definition}
\newtheorem{rema}[equation]{Remark}

 \numberwithin{equation}{section}

\begin{document}
\title[Automorphism groups of Steiner triple systems]
{Automorphism groups of Steiner triple systems
 }
  
       \author{Jean Doyen}
        \address{Universit\'e Libre de Bruxelles, Bruxelles 1050, Belgium}
     \email{jdoyen@ulb.ac.be}
        \author{William  M. Kantor}
       \address{U. of Oregon,
       Eugene, OR 97403
        \      and  \
       Northeastern U., Boston, MA 02115}
       \email{kantor@uoregon.edu}
  
  \begin{abstract} 
    \label{corollary}
If $G$ is a finite group  then there is an integer $M_G$ such that$,$ 
for  $u\ge M_G$ and 
$u\equiv 1$ or $3$ {\rm(mod~6),}   there is a Steiner triple system $U$ on $u$ points
  for which   ${\rm Aut}\hspace{.5pt} U \hspace{.5pt}\cong G.$ \
If $V$ is a Steiner triple system then there is an integer $N_V$ 
such that$,$  for  $u\ge N_V$ and 
$u\equiv 1$ or $3$ $($mod~$6),$   there is a Steiner triple system $U$ on $u$ points
having $V$ as an  ${\rm Aut}\hspace{.5pt} U$-invariant   subsystem  
such that   ${\rm Aut}\hspace{.5pt} U\cong{\rm Aut} V $
and   ${\rm Aut}\hspace{.5pt} U  $  induces  ${\rm Aut} \hspace{.5pt}V$ on $V$.
 \end{abstract}

\maketitle

  \section{Introduction}
  \label{Introduction}
  Mendelsohn \cite{Me} proved that any finite group $G$ is isomorphic to the   automorphism
  group of some Steiner triple system.  In his proof he   modified
   the Steiner triple system of points and  lines of a projective space $PG(n,2)$,
   producing a system having $2^{n+1}-1$ points for some $n$. 
  This leads to the natural question: 
  what restrictions are  there  on the number of points of a Steiner triple system
  $U$ such that ${\rm Aut}\hspace{.5pt} U\cong G$?  
  In order to admit $G$ as a group of automorphisms,   $U$  cannot be too small:  
\begin{theorem}
  \label{corollary}
If $G$ is a finite group  then there is an integer $M_G$ such that$,$ 
for  $u\ge M_G$ and 
$u\equiv 1$ or $3$ {\rm(mod~6),}   there is a Steiner triple system $U$ on $u$ points
  for which   ${\rm Aut}\hspace{.5pt} U  \cong G.$ 
\end{theorem}

As with most theorems of this sort, the proof does not distinguish between cyclic groups
and simple  groups.  It is  known that $M_G=15$ when $G=1$ \cite{LR}.
Our arguments cannot deal with such small Steiner triple systems.

The preceding theorem  is an~immediate consequence of  \cite{Me}  and the following   more general result,
which is the main theorem of this paper:

\begin{theorem}
  \label{main}
If $V_{\scriptscriptstyle\star} $
is a Steiner triple system then there is an integer 
$N_{V_{\star}}$\! such that$,$ 
for  $u\ge N_{V_{\star}}\!$ and 
$u\equiv 1$ or $3$ {\rm(mod~6),}   there is a Steiner triple system $U$ on $u$ points having ${V_{\scriptscriptstyle\star}}$ as an  
${\rm Aut}\hspace{.5pt} U$-invariant   subsystem  
such that   ${\rm Aut}\hspace{.5pt} U\cong{\rm Aut} {V_{\scriptscriptstyle\star}} $
and   ${\rm Aut}\hspace{.5pt} U  $  induces  ${\rm Aut} \hspace{.5pt}{V_{\scriptscriptstyle\star}}$ on ${V_{\scriptscriptstyle\star}}$.
\end{theorem}

Cameron \cite{Ca} considered a similar question.  He proved that,
if ${V}$ is a Steiner triple system of order $v$ (i.\,e., having $v$ points),
and if $u>6v^2$ with 
$u\equiv 1$ or $3$ $($mod~$6),$  then
 there is~a Steiner triple system $U$ of order $u$ in which ${V}$ can be  embedded 
in such a way~that every automorphism of ${V}$  can be extended to $U$. 
His proof and ours  use a familiar and wonderful  construction of Moore
\cite{Mo} 
from the not-so-distant past that combines three Steiner triple systems to produce a fourth. 

The proof of Theorem~\ref{main}  first enlarges $V_{\scriptscriptstyle\star}$,  without changing its automorphism group,  as part of a process to obtain 
a Steiner triple system $U$
 having a  rich   geometry of  $PG(m,2)$  subsystems for various $m$
 (cf.~Remark~\ref{rich}).
 This   process  involves  Lemma~\ref{direct product}  and 
 Proposition~\ref{VxD} (using \cite{Ka}),  and  leads to our key tool: Proposition~\ref{maximal}.
  The latter makes it straightforward in  Proposition~\ref{Aut A}
  to determine the automorphism group of the 
  Steiner triple system $U$ we construct.

The  ugly bookkeeping parts of the proof  
(in Section~\ref{Properties of $X,$ $Y$ and $V$} and especially 
 in Section~\ref{Existence of $y$ and $x$})
 ensure  that we obtain all large $u$.   
 Remark~\ref{bounds}  contains poor bounds for $M_G$ and~$N_{V_{\star}}$,
while     Remark~\ref{comparison}  comments on a difference between 
 the bookkeeping approaches in \cite{Ca} and here. 
    
 Results such as Theorem~\ref{corollary} are usually based on the action of $G$ 
  having many regular point-orbits.  This is very much not the situation for 
 Theorem~\ref{main}: for our $U$  the size of every point-orbit of   
 ${\rm Aut}\hspace{.5pt} {U}$   is
  1  or the size of a point-orbit of 
  ${\rm Aut} {V_{\scriptscriptstyle\star}}$ on  the  original 
  Steiner triple system
  ${V_{\scriptscriptstyle\star}}$.
 
  There is also  a result  in \cite{Ca}  concerning a {\em partial Steiner triple system} $V$ 
  (a set of points, together with some triples of points, such that any two points are in at most one triple), and a partial Steiner triple   $\,U$ having $V$ as a 
{\em  subsystem}   (so the points of $V$ are  among the points of $ U$, and 
the triples of $V$ are  precisely those triples of  $U$ that are 
contained in $V$).  
 It is shown   in \cite{Ca}  that  there is a function $g$ such that,
   if ${V}$ is a partial Steiner triple system of order $v,$ and if $u>g(v)$ with 
$u\equiv 1$ or $3$ $($mod~$6),$  then
 there is~a Steiner triple system $U$ of order $u$ of which ${V}$ is a subsystem 
  such  that every automorphism of ${V}$  can be extended to $U$. 
In Section~\ref{Proof PSTS} we will use    Theorem~\ref{main} to prove the following stronger  result (along with variations):

\begin{theorem}
  \label{PSTS}
If $V  $  is a partial Steiner triple system then there is an integer 
$N'_{V }$  such that$,$ 
for  $u\ge N'_{V }$  and 
$u\equiv 1$ or $3$ {\rm(mod~6),}    there is a Steiner triple system $U$ on $u$ points
having $ V $ as an  ${\rm Aut}\hspace{.5pt} U$-invariant   subsystem   
such that ${\rm Aut} \hspace{.5pt}U\cong{\rm Aut} V $ 
and
${\rm Aut} \hspace{.5pt}U  $  induces  ${\rm Aut} V$ on $V$.
 \end{theorem}

  \section{Background}
  \subsection{Moore's $XYV$}
  \label{XYV}
  We will use a 125 year old construction due to Moore
   \cite[p.~276]{Mo}\footnote{Moore used it to
  produce  two nonisomorphic STSs of any admissible order
   $>13$.     Unfortunately, his method for  proving nonisomorphism
    \cite[pp.~279-281]{Mo}  has a significant gap.}.
  (This  construction is  in many sources, such as \cite[p.~235]{LW} and  \cite{Ca}.)   
  
  Let $X\subset Y$ and $V$ be three STSs (i.\,e., Steiner triple systems), and label $Y-X$  
  in any way by the elements of a cyclic group  $A$ of order $|Y|-|X|$.
   (We always use $|Y|  $ to denote the order of an STS  $\,Y$.)
  Then $U:=X\cup (V\times A)$ is (the set of points of) an STS, with triples
  \begin{itemize}
  \item[(M1)] 
  those of $X$, 
  \item[(M2)] 
$(v,a_1), (v,a_2)$ and 
 $\begin{cases} 
 x&$\!\!$ \mbox{if $a_1,a_2,x$ is a triple in $Y$, $a_i\in Y-X$, $x\in X$}
 \\
 (v,a_3)& $\!\!$\mbox{if $a_1,a_2,a_3$ is a triple in $Y$, $a_i\in Y-X$, }
 \end{cases}$  
 
 \noindent
 and
  \item[(M3)] $(v_1,a_1), (v_2,a_2 ), (v_3,a_3)$ if $v_1,v_2,v_3$ is a triple in $V$ and
  $a_1a_2a_3=1$.
   \end{itemize}

The fact that $A$ is cyclic, not just abelian, is used in several places, most significantly in  Lemma~\ref{labeling}   and Proposition~\ref{Aut A}.
        
 Clearly  $|U|=|X|+|V|(|Y|-|X|)$. 
 
  \subsection {Enlarging $Y$} 
    \label {$2Y+1$ construction.} 
    \label{2Y+1}
    The STSs $X$ and $Y$ in the preceding section have unknown structure.  
    While this does not matter for $X$, 
    we will use  elementary constructions to enlarge $Y$  in order to give it  significant geometric structure (Lemmas~\ref{PG(2,2)-pointed} and~\ref{direct product}).
        
  Given an STS   $\, Y_0$ there is  a standard construction for an STS $\, 2Y_0+1$
    on  $2|Y_0|+1$ points, labeled using
  $Y_0\dot\cup Y_0'\dot\cup *_0 $  for a ``distinguished'' 
  new point $*_0$  and
    a  bijection $y \mapsto y'$ sending   $Y_0 \to Y ' _0$, and  
  triples of the~form
 $$\mbox{ \ \  \ \ 
  $abc$ in $Y_0$
   \ \  \ \ 
  $*_0aa'$ etc.
   \ \  \ \ 
  $a'b'c$ etc.}
  $$
  
  \noindent
  Here $|2Y_0+1|\equiv 3$~(mod~4) and $Y_0$ is a subsystem of  $2Y_0+1$.
  If $Y_0$ is a hyperplane of a projective space $P=PG(n,2)$ 
then $P\cong 2Y_0+1$.   Thus, if $Y_0$ is a projective space then so is $2Y_0+1$.


  \defi\rm An STS is~{\em $PG(2,2)$-pointed}  with respect to a   point $p$
if any two triples containing $p$ generate a $PG(2,2)$ subsystem.

 An STS   with more than seven points is {\em $PG(3,2)$-$2$-pointed}  with respect to two   points 
if any four points including these two   generate a $PG(k,2)$ subsystem for $k=2 $ or 3.
  
    An STS   is \emph{$PG(3,2)$-paired}  if any two points are in a
    $PG(3,2)$ subsystem.   {\em   This is the key geometric property}
 {\em   needed in the proof of} Proposition~\ref{maximal}(i). 
  
 \begin{lemma}
   \label{PG(2,2)-pointed}
   If $Y_0$ is an STS with more than one point$,$ then
 \begin{itemize}
 \item[(i)]   
 $2Y_0+1$  is $PG(2,2)$-pointed with respect to the distinguished point $*_0,$   
 \item[(ii)]    $Y_1:=2(2Y_0 + 1) + 1$  is  
$PG(3,2)$-$2$-pointed   {$($with respect to some pair of its points$),$}  and
  
   \item[(iii)]  $Y_1$ is $PG(3,2)$-paired of order      
  $ {4|Y_0 |+3}\equiv 7$~{\rm (mod~$8$).}

  \end{itemize}
 \end{lemma}
 \proof
 (i)  
 Two triples $*_0aa'$, $*_0bb'$ of $ 2Y_0 + 1 $ containing $*_0$ generate a $PG(2,2)$ subsystem with triples
 $$\mbox{  \ \   \ \ 
  $abc$
   \ \   \ \  \ \ 
  $*_0aa'$, $*_0bb'$, $*_0cc'$
   \ \  \ \  \ \ 
   $a'b'c$,  $a'bc'$,  $ab'c'$. } 
   $$ 

(ii) The   pair  $\{ *_0,*_1\}$ has the required property, 
where  
$*_1$ is the new point used to produce $ Y_1= {2 (2Y_0 +1)}+1  $ from $2Y_0 +1$. 
For, if $*_0,*_1,a,b$ are 
four points of $Y_1$  then $*_1,a $  and  $*_1, b$  
are in triples meeting $2Y_0 +1 $ at points $\bar a $ and $\bar b$, respectively.
Then $*_0,\bar a,\bar b$ are in a $PG(2,2)$ subsystem $Z$ of $2Y_0+1$,
and $2Z+1$ is a $PG(3,2)$ subsystem of $ Y_1 $  containing $*_0,*_1,a,b$. 

(iii)  This is immediate by (ii).
  \qed
           
       \medskip
{\em Admissible integers} are those  $\equiv1$ or 3 (mod 6); these are precisely the possible  orders of STSs.

   We use  the standard {\em direct product} 
   \mbox{$A{\bm\times}B$} 
   of STSs
$A$ and  $B$:     the STS for which $A{\times} B$ is the set of points and
whose triples    have the form 
$\{(a,b_1),(a,b_2),(a,b_3 ) \}$,
$\{ (a_1,b), (a_2,b), (a_3,b)\}$
  or 
$\{(a_1,b_1),(a_2,b_2),(a_3,b_3)\}$
for $ a\in A,  $  $b\in B$, and ordered triples $(a_1,a_2,a_3)$ from $A$ and 
$(b_1,b_2,b_3)$ from $B$.
Clearly, if $a\in A$ and $b\in B$ then $a{\bm\times}B$  and $A{\bm\times}b$
are subsystems of $A{\bm\times}B$.
  
  \begin{lemma}
   \label{direct product}~
     \begin{itemize}
   \item [\rm(i)]     
If $A$ and $B$ are $PG(3,2)$-paired  STSs then so is 
 \mbox{$A{\bm\times}B$}.    
     \item [\rm(ii)]  An  integer is  the  order of a $PG(3,2)$-paired  STS
   if it is admissible$,$    at least  $15$   and  ${\equiv 7}$~{\rm(mod~8),}    
   or  if it is  a product of   integers 
    each of which  behaves that way.
   \end{itemize}
   \end{lemma}
   
   \proof (i)
Given two points  $(a_1,b_1)$ and $(a_2,b_2)$  of  \mbox{$A{\bm\times}B$}  there are  
$PG(3,2)$ subsystems   $A_0$ of $A$ containing $a_1,a_2$ and
 $B_0$ of $B$ containing $b_1,b_2$.

If $a_1\ne a_2$ and $b_1\ne b_2$   there is an isomorphism $\phi\colon \!A_0\to B_0$
sending $a_i\mapsto b_i$  for  $i=1,2$.~Then 
$\{(a ,a ^\phi ) \! \mid \!a\in A_0  \}\subset A_0{\bm\times} B_0$
 is a $PG(3,2)$ subsystem of  \mbox{$A{\bm\times}B$} containing the given~points.

If $a_1= a_2$ then  
 $a_1 {\bm\times} B_0$ is a $PG(3,2)$ subsystem
 containing $(a_1,b_1)$ and $(a_2,b_2)$.  If 
$b_1= b_2$ then  $A_0  {\bm\times}  b_1$  is a $PG(3,2)$ subsystem
 containing $(a_1,b_1)$ and $(a_2,b_2)$.

\smallskip
(ii)  Any admissible integer $8n+7 \ge 15$ can be written 
$2(2(2n+1)+1)+1$ with $2n+1>1$ and
admissible.  By Lemma~\ref{PG(2,2)-pointed}(iii)
there is a $PG(3,2)$-paired  STS
of order $8n+7$.  Now use (i).
\qed 
        
    

    \subsection{Enlarging ${V_{\scriptscriptstyle\star}}$}
  \label{Increasing $V$.}
The STS   $\,  {V_{\scriptscriptstyle\star}}$ in Theorem~\ref{main}
 has   unknown structure.
As  we did with $Y_0$ in Section~\ref{2Y+1}, 
we will enlarge  $\,  {V_{\scriptscriptstyle\star}}$
 to STSs having   significant geometric structure  (Corollary~\ref{=3}).
 Since our arguments are based on finite geometry, 
 we briefly describe   how ``close'' each  
 STS  in \cite{Ka} is to a projective space.   
 
 Let $V$ be a  vector space  over  $F={\bf F}_{16}$ with basis $v_1,\dots,v_n.$
This produces a vector space and  a projective space  over ${\bf F}_2$.  Let $F^*=\langle \theta\rangle$.
 Suitably modify the  $PG(3,2)$  subspaces determined by the ${\bf F}_2$-spaces 
 $Fv_i$ and $F(v_i+\theta v_j)$ for all $i,j,$ in order to obtain an STS $D$.   Only a tiny portion
 of the underlying  $PG({4n-1},2)$ is altered:  every subspace of  that projective space
  meeting $U:=\bigcup_{i,j} \big (Fv_i\cup {F(v_i+\theta v_j)} \big)$  at most once is a subsystem of $D$.  As in \cite[Theorem 1.1 and Sec.~7(1c)]{Ka}, this provides the flexibility needed  for the following

\rema 
\label{D}
If $n\ge 6$  then    there is an STS $\,D$ such that  
  \begin{itemize}
\item [(i)]${\rm Aut} D=1,$
\item [(ii)]  $D$ has $16^n-1$ points$,$ and 

\item [(iii)]  given points $a,b$ of $D$ there is a point $c$ such that each   pair 
$\{a,c\}$ and $\{b,c\}$ is in some $PG(n-1,2)$~subsystem of $D$.
 
\end{itemize}\rm

See \cite{Ka} for the modifications of the above $PG(3,2)$ subsystems needed to obtain (i).  
The   subsystems in (iii) are crucial for  our proof of Theorem~\ref{main};
we will   obtain  them  less tediously than   similar ones in \cite{Ka}.   
Let $c\notin \bigcup_{u\in U}   
( \langle a,u \rangle _{{\bf F}_2}     \cup \langle b,u \rangle _{{\bf F}_2}  )$ and consider the pair $\{a, c \}$.  
For $2\le j<n$  inductively increase a $j$-dimensional ${\bf F}_2$-subspace $J\supset\{a, c \}$
of $V$ with  $J\cap U\subseteq  \langle a \rangle _{{\bf F}_2} $  to
a $j+1$-dimensional ${\bf F}_2$-subspace $J'\supset J$ with $J'\cap U=J\cap U$,
noting that  $|\bigcup _ { u\in U }\langle J,u \rangle _{{\bf F}_2}|
< 2^{ j+1}(n+n(n-1)) 16 <|V|$.

\smallskip\smallskip

{\bf  Choose $n\ge 6$ so that $ 2^{n }-1>|{V_{\scriptscriptstyle\star}}| \ge 2^{n-6} $.}		 Then  \emph{the image of  any map~from $PG(n-1,2)$ into ${V_{\scriptscriptstyle\star}}$ sending every collinear triple to a triple or a point must have   size} 1.
(Otherwise, since the map cannot be 1-1,
 restrict to a plane mapping onto a triple in order to obtain a contradiction:
the preimages of the points of the triple would have to be pairwise disjoint,
 cover the plane, and be such that the line through two   points of a preimage
 is contained in the preimage.)
 
 \begin{prop}
\label{VxD} 
Let $D$ be as in the preceding {\rm Remark}$,$ and let $d\in D $.
   \begin{itemize}
    \item[(i)]
     Every point of 
   ${V_{\scriptscriptstyle\star}}{\bm\times} D$ is in a $PG(n-1,2)$ subsystem of 
     ${V_{\scriptscriptstyle\star}}{\bm\times} D,$
     
   \item[(ii)]
 ${V_{\scriptscriptstyle\star}}{\bm\times} d $ is isomorphic to ${V_{\scriptscriptstyle\star}},$   
 
   \item[(iii)]
   ${V_{\scriptscriptstyle\star}}{\bm\times} d $ is an ${\rm Aut}({V_{\scriptscriptstyle\star}}{\bm\times} D)$-invariant subsystem of 
 ${V_{\scriptscriptstyle\star}}{\bm\times} D$    
 on which 
${\rm Aut}({V_{\scriptscriptstyle\star}}{\bm\times}D)$  induces  
$  {\rm Aut} 
  ({V_{\scriptscriptstyle\star}}{\bm\times} d),$
  and
  $ {\rm Aut} ({V_{\scriptscriptstyle\star}}{\bm\times} D)  \cong {\rm Aut} 
  ({V_{\scriptscriptstyle\star}}{\bm\times} d)   
   \cong {\rm Aut} \hspace{.5pt} {V_{\scriptscriptstyle\star}}  ,$
   and
    \item[(iv)]  $|{V_{\scriptscriptstyle\star}}{\bm\times} D|\equiv \pm3$
    {\rm (mod~12)}  and  
$16^5<|{V_{\scriptscriptstyle\star}}{\bm\times} D|<2^{ 24} 
| {V_{\scriptscriptstyle\star}}|^5$. 
\end{itemize}
 \end{prop}
 
 \proof 
For (i), if $v\in {V_{\scriptscriptstyle\star}}$  then  $v{\bm\times}D $  is isomorphic to $D$, so use Remark~\ref{D}(iii).  Statement (ii)  is obvious.  
For (iii) we need to determine  ${\rm Aut}({V_{\scriptscriptstyle\star}}{\bm\times}D)$.
 
  Every  $PG(n-1,2)$ subsystem  in ${V_{\scriptscriptstyle\star}}{\bm\times}D$ 
  is a set of ordered pairs that projects onto  subsystems of
  ${V_{\scriptscriptstyle\star}}$ and $D$,  and hence induces a map from $PG(n-1,2)$ to $V_{\scriptscriptstyle\star}$
  sending every collinear triple  to a point or triple.  As noted above,
  since $2^n-1> |{V_{\scriptscriptstyle\star}}|  $ the image of that map is a point  of   ${V_{\scriptscriptstyle\star}}$.  
  Thus, {\em every $PG(n-1,2)$ subsystem  in ${V_{\scriptscriptstyle\star}}{\bm\times} D$ 
  lies in some subsystem $v{\bm\times} D,$ $ v\in {V_{\scriptscriptstyle\star}}$}.

 Consider the graph whose vertices are the  $PG(n-1,2)$ subsystems of 
 ${V_{\scriptscriptstyle\star}}{\bm\times} D$, with two  such   subsystems
  joined if and only if they meet.   We have just seen that every such   subsystem is contained in some subsystem
  $v{\bm\times} D$, $v\in V$.
Since no two subsystems $v{\bm\times} D$ meet, every connected component $C$ of our graph 
lies in  some  subsystem $v{\bm\times} D;$ 
 by Remark~\ref{D}(iii),  $C$
 is the set of $PG(n-1,2)$ subsystems of   $v{\bm\times} D$  
 and   generates $v{\bm\times} D$.
Since   ${{\rm Aut}({V_{\scriptscriptstyle\star}}{\bm\times}D) }$ 
permutes the connected components of our graph it also 
 permutes  the  subsystems $v{\bm\times} D$.
 
  Let $h\in {\rm Aut}({V_{\scriptscriptstyle\star}}{\bm\times} D)$.
 If $v\in V_{\scriptscriptstyle\star}$ then 
 $(v{\bm\times}D)^h= v^{h'}{\bm\times}D $ for some 
 $v^{h'}\in V_{\scriptscriptstyle\star}$, where 
 $h'\colon V_{\scriptscriptstyle\star}\to V_{\scriptscriptstyle\star}$ is bijective.  
 We claim that the map
 $h'\colon \!{V_{\scriptscriptstyle\star}}\to {V_{\scriptscriptstyle\star}}$
 {\em is in} ${\rm Aut}V_{\scriptscriptstyle\star}$.~~For, 
 if $v_1,v_2,v_3$ is a triple of $V_{\scriptscriptstyle\star}$ 
 and $d\in D$ then 
 $(v_1,d),(v_2,d),(v_3,d)$ is a triple of 
 ${V_{\scriptscriptstyle\star}}{\bm\times} D$.
 Then so is $(v_1,d)^h,(v_2,d)^h,(v_3,d)^h$;
 this is $(v_1^{h'},d_1),(v_2^{h'},d_2),(v_3^{h'},d_3)$  
 with $d_i\in D$,  so $v_1^{h'},v_2^{h'},$ $ v_3^{h'}$
 is a triple of  $V_{\scriptscriptstyle\star}$, as claimed.
 
 Now $h' \in {\rm Aut}\hspace{.5pt} V_{\scriptscriptstyle\star}$ induces
 $h^\bullet \in   {{\rm Aut}({V_{\scriptscriptstyle\star}}{\bm\times}d) }$ 
  sending   
 $(v,d)\mapsto (v^{h'},d)$ for ${v\in {V_{\scriptscriptstyle\star}}},$ ${d\in D}$.
 We thus have two automorphisms $h$   and $h^\bullet$ of 
$ {V_{\scriptscriptstyle\star}}{\bm\times}d $ sending each 
 $v{\bm\times} D $ to  $v^{h'}{\bm\times} D $.
 Then $ h^\bullet  h^{-1}$ sends each subsystem $ v{\bm\times} D $
 to itself,  induces
    an automorphism of each such subsystem, and hence is 1   by (ii) and Remark~\ref{D}(i).
    Thus, $ h=h^\bullet $ sends each subsystem 
    $ {V_{\scriptscriptstyle\star}}{\bm\times} d $
 to itself, and hence so does ${\rm Aut}({V_{\scriptscriptstyle\star}}{\bm\times}D)$.
 This proves that ${\rm Aut}({V_{\scriptscriptstyle\star}}{\bm\times}D)=
( {\rm Aut}{V_{\scriptscriptstyle\star}})\times 1_D$.

(iv)   $|V_{\scriptscriptstyle\star}{\bm\times}D|
 = (16^n-1)|{V_{\scriptscriptstyle\star}}|\equiv 
(16-1) |{V_{\scriptscriptstyle\star}}|\equiv
 3|{V_{\scriptscriptstyle\star}}|\equiv
\pm3$~{\rm (mod~12)}
 and  $16^5<{(16^n-1)}|{V_{\scriptscriptstyle\star}}|<
 (2^n)^4| {V_{\scriptscriptstyle\star}}  | \le
(2^{ 6}  |  {V_{\scriptscriptstyle\star}} |)^4 | {V_{\scriptscriptstyle\star}}  |$
since $| {V_{\scriptscriptstyle\star}}  |\ge 2^{n-6}$.
   \qed 
 
 \begin{corollary}
 \label{=3}
 For   sufficiently large integers $n$ and $m$ there are STSs 
 $V_{\scriptscriptstyle\star}{\bm\times}D$
 and  $V_{\scriptscriptstyle\star}{\bm\times}D{\bm\times}D'$
 such that one of them$,$ $V,$ has the following properties$:$
  \begin{itemize}
   \item[(i)]  $|V|\equiv 3 $   {\rm (mod~12),}  
   $|V|>16^5,$
   and  either
  $|V|=|V_{\scriptscriptstyle\star}||D|=|V_{\scriptscriptstyle\star}|
   (16^n-1) $   or 
   $|V|=|V_{\scriptscriptstyle\star}||D||D'|=|V_{\scriptscriptstyle\star}|
   (16^n-1)(16^m-1) ,$
    \item[(ii)]
     Every point of 
   $V$ is in a $PG(3,2)$ subsystem of  $V,$
     and
   \item[(iii)]
$ V$ has an ${\rm Aut}\hspace{.5pt} V$-invariant subsystem  
${V_{\scriptscriptstyle\star}}' 
\cong{V_{\scriptscriptstyle\star}} $    
 on which 
${\rm Aut}\hspace{.5pt} V$  
induces~${\rm Aut}\hspace{.5pt} {V_{\scriptscriptstyle\star}}' ,$
  and
  $ {\rm Aut} V\cong {\rm Aut}  {V_{\scriptscriptstyle\star}}'
  \cong {\rm Aut}  {V_{\scriptscriptstyle\star}}.$
  \end{itemize}
  \end{corollary}
  
  \proof
  Apply the proposition to ${V_{\scriptscriptstyle\star}}{\bm\times} D$
  in place of ${V_{\scriptscriptstyle\star}}$,
  using in place of $D$ an STS $D'$ of order $16^m-1$ 
  where $2^m-1>|{V_{\scriptscriptstyle\star}}{\bm\times} D| \ge 2^{m-6}$.
  Since $3  |V_{\scriptscriptstyle\star}|\equiv\pm3$ (mod 12) 
  either $ |V_{\scriptscriptstyle\star}|(16^n-1)\equiv3$ (mod 12)   or 
   $ |V_{\scriptscriptstyle\star}|(16^n-1)(16^m-1) \equiv3$ (mod 12).   \qed

     \section {Proof  of  Theorem~\ref{main}} 
     \label {Proof} 
     The proof proceeds in three stages:  construction of an STS $U$ 
     (Section~\ref{Properties of $X,$ $Y$ and $V$}),  determining that we have obtained all sufficiently large 
     admissible integers as the order of some such $U$   (Section~\ref{Existence of $y$ and $x$}), 
     and using geometry to  determine {\rm Aut\hspace{.5pt}}$U$  (Sections~\ref{Finding critical subsystems} and \ref{The occurrence of  $A$.}).
     
  \subsection {Preliminaries and notation} 
  \label{Preliminaries} 
  We  will describe   STSs~$X,$ $Y$ and  $V$ that will be used to construct our STS~$U$  via  Section~\ref{XYV}.
  
 \subsubsection {Properties of $X,$ $Y$ and $V$}
 \label {Properties of $X,$ $Y$ and $V$}
  We begin with notation and properties of these STSs.
  We note that properties  (P1(b)) and (P3(d))
  will be essential  (in Proposition~\ref{maximal})
     for studying  subsystems of $U$ isomorphic to $V$ or $Y$.
  
  Let $V_*$ be as in Theorem~\ref{main}.
  Property (P1) concerns an STS $\,V$ that will replace $V_*$  in our arguments
  and has a rich geometric structure.   
      \begin{itemize}
      \item [(P1)]    {\em$V$ and $v$.}
        \begin{itemize}
        \item     [(a)] 
Use   Corollary~\ref{=3} to obtain an STS $\,V$     on  $v$ points, where

\smallskip
\qquad\qquad $v\equiv 3 $   {\rm (mod~12)}  and $ v>16^5$, 
\smallskip

having  (a copy of) ${V_{\scriptscriptstyle\star}}$ as an
${\rm Aut} V $-invariant  subsystem on which 
    ${\rm Aut} \hspace{.5pt} V $ induces ${\rm Aut}\hspace{.5pt}  {V_{\scriptscriptstyle\star}}$ 
    and such that  ${\rm Aut} \hspace{.5pt} V \cong {\rm Aut} {V_{\scriptscriptstyle\star}}.$
(Here $V $ is not uniquely determined: it depends on   choices  for $D,$ $n$  and possibly $D'$ and $m$ that are   made just once in  the proof of Corollary~\ref{=3}.)
    \item     [(b)] 
 Every point of $V$ is in a  $PG(3,2)$ subsystem of $V$
  (Corollary~\ref{=3}(ii)).   
    \end{itemize}
      \end{itemize}      
 In (P2) and (P4) we will  introduce further admissible integers $y_1$,
        $y_2$,  $y$, $x$  and $u$;   Lemma~\ref{Existence2} concerns the  existence of 
         integers satisfying the conditions stated in (P1),  (P2) and (P4).
        
        
    \begin{itemize}
   \item [(P2)]
  {\em$\delta, $  $x$, $y_1,$ $ y_2$ and $y$.}
    \begin{itemize}
  \item     [(a)]   Let $\delta=\pm1$.   (The admissible integer $u$ in 
  Theorem~\ref{main} will later be related to~$\delta$ via  the requirement
    $u  \equiv\delta $~(mod~4).)
  
  We will use an admissible integer $y \equiv\delta $~(mod~4):
  
\ \ \  let  $y_1\equiv 15$~(mod~24), which is admissible   and $  \equiv 7$~(mod~8);
    
\ \ \ if $\delta=-1$  let $y:=y_1$; and
        
\ \ \ if $\delta= 1$  let $y:=y_1y_2$ where $y_2\equiv15$~(mod~24).

 \item     [(b)]  Let $x$ be admissible.
 
 \item     [(c)] Assume that  $y_1\ge \sqrt y\ge 8x+7$. 
 
  \item     [(d)]  Assume that   $y>x+  6 v$.
             \end{itemize}
                \end{itemize}  
                  
        \noindent
        In (P3) and (P4) we provide a recipe that uses the  integers in (P2) to obtain 
 auxiliary STSs together with an STS $\,U$ that  behaves as required in Theorem~\ref{main}.
        
  \begin{itemize}
   \item [(P3)]    {\em$X,$  $Y_1,$ $Y_2$ and $Y$.}
 \begin{itemize}
        \item     [(a)]  
 Write $y_1=4y_0+3$, so that    $y_0$ is admissible and
$y_0\ge2x+1$ (by (P2(a))  and  (P2(c))).  Then
  \cite{DW} provides an  STS   $\,Y_0$ of order  $y_0$ 
  containing a subsystem  $X$ of order $x$.     
    \item     [(b)] 
  Let $Y_1:=2(2Y_0+1)+1$, so  $ |Y_1|=2(2|Y_0|+1)+1=y_1$  (by (P3(a)))
  and $Y_1$ is $PG(3 ,2)$-paired (by Lemma~\ref{PG(2,2)-pointed}(iii)).
  
 If  $\delta=1$  let   $Y_2$ be a $PG(3 ,2)$-paired STS of order $y_2$ 
(Lemma \ref{direct product}(ii)).
 \item     [(c)]    
If $\delta=-1$  let $Y:=Y_1$,   of order $y$.
        
If $\delta=1$  let $Y:=Y_1{\bm\times} Y_2$,   of order $y$.

\item     [(d)]  Any two points of  $ Y$ are in a  $PG(3,2)$ subsystem
of  $ Y$ (Lemma~\ref{direct product}).
   \end{itemize}  
       \end{itemize}

      \begin{itemize}
 \item [(P4)]
    {\em$U,$ $A$ and $u$.}
 
 \noindent
  Let $  A:= Y-X$  be as in Section~\ref{XYV}, so the STS    $U:=X\cup (V\times  A)$
 has order~$u:=x+v(y-x)$.   In  Lemma~\ref{Existence2} we will see that 
 all sufficiently large admissible integers  arise here  as $u$  
 for the choice of $v$   in (P1)  and  for suitable $x,y,\delta$ in  (P2).
 
 \noindent
As noted in   \cite[p.~469]{Ca},  each $g\in {\rm Aut}  V$ acts 
     as an automorphism of  $U$ via%
     \vspace{-2pt}
$$
 \mbox{$g=1$ on $X$ and $(p,q)^g=(p^g,q)$ for $p\in V,$ $ q\in  Y-X$.}\vspace{-2pt}
 $$
 This produces a subgroup of 
 ${\rm Aut}\hspace{.5pt}   U$ isomorphic to ${\rm Aut}  \hspace{.5pt} V$
 and inducing ${\rm Aut}  \hspace{.5pt} (V\times 1)$  on the subsystem 
 $V\times 1$ of $U$.

 \smallskip
 
 \item [(P5)]
  The cyclic group $ A$ has even order; let $-1$ denote its involution.
  For $a\in  A$ let $-a:=(-1)a$.
   
Let $ A_6:=\{a\in   A\mid a^6=1\}$. 
\smallskip

  \item [(P6)]
  {\em Labeling $  Y-X$.}
    We  assume that $Y-X$  behaves  as in Lemma~\ref{labeling} 
  below (which depends only on Section~\ref{XYV} and the fact that 
  $|Y-X|$ is not tiny).
  \end{itemize}
 
\begin {rema}
\label{rich} 
\rm Section~\ref{Introduction}  mentions ``\em\dots a Steiner triple system $U$
 having a  rich   geometry of  $PG(m,2)$  subsystems\dots\rm''.   
  This refers to the  geometry  inherited by $U$ from   (P1(b)) and (P3(d)), which    involves   far more structure  than the fact that $U$ is generated by its $PG(3,2)$  subsystems.
\end {rema}

 \subsubsection {Existence of $y$ and $x$}
 \label {Existence of $y$ and $x$}
We first rephrase and slightly strengthen the numerical requirements in 
(P1),   (P2) and (P4):
 
 \begin{lemma}
   \label{Existence1}
 Assume that
 $u,v,\delta , y_1,y_2$  and $y$  are  integers that   behave as follows$:$
\begin{enumerate}
\item [\rm(i)] 
$u$  is admissible with $u \equiv \delta {~\rm(mod~4)}$ for  $\delta = \pm1,$
and $v\equiv 3$ {\rm(mod 12)}$,$
\item [\rm(ii)]    $u> 800^2v^7$   and  $v>16^5 , $  
\item [\rm(iii)] 
$v-1$ is a factor of $u-y , $    
\item [\rm(iv)]
 $y=y_1y_2 \equiv \delta {~\rm(mod~4)}$ with $y_1 \equiv 15${~\rm(mod~24)}$ ,$   where

\noindent
 if  $\delta=-1$ then $y_2=1,$ while 

\noindent
 if  $\delta=1  $   then $y_ 2 \equiv 15${~\rm(mod~24)} and  $y_1\ge y_2 $
  $($so  $y_1\ge \sqrt y\hspace{1pt}) ,$
and
 \item [\rm(v)] 
 ${u/v}< y <{u/v} +\frac{1}{8}\{(v-1)/v\}\sqrt {{u/v}}\, -1 <u$.
 
\end{enumerate}
Then $u,v,\delta, y_1,y_2,y$ and
$  x :=(vy-u)/(v-1)  =y- (u-y)/(v-1) $  are integers  that satisfy all of the conditions in 
{\rm (P1),}   {\rm (P2)}  and {\rm (P4).}
\end{lemma}

{\em\noindent Remark.} 
We have   $y_1 \equiv 7$~{\rm(mod~8)}, 
and   $y_2 \equiv 7$~{\rm(mod~8)}  if   $y_2\ne 1$.  
However, we do not need information about either $u$ or $v$~{\rm(mod~8)}.   
What we need are
$u $ and~$y$~{\rm(mod~4)} in (i) and (iv) in order to have 
$  u - y \equiv0$~(mod~4);   this  and $v-1\equiv  2$  (mod 4)  imply that  
  $ x=y-(u-y)/(v-1) \equiv  y-0\equiv1$~(mod~2).%

\proof
Note that  $u,v,y_1,y_2$  and $y$ are admissible.
By (v), $vy-u> 0$  and $u- y>0$, so $0<x <y$.

\smallskip\smallskip
(P1):  See (i) and (ii).

\smallskip\smallskip
(P2(a)):  This is in (iv).
\smallskip\smallskip

(P2(b)):   Since $v\equiv 0$ (mod 3) we have 
$x \equiv(0-u)/(0-1)=u\equiv 0$ or $1$~(mod~3).  
We have already noted that $x$ is odd, so it is admissible.
\smallskip\smallskip

(P2(c))  and (P2(d)):    By (v),  $(vy-u)/v 
 < \frac{1}{8}\{(v-1)/v\}\sqrt {{u/v}}\, - \frac{7}{8}\{(v-1)/v\}$.~~Then  
\begin{equation*}
  \begin{array}{lllllll}
  x  &\hspace{-6pt}=\hspace{-6pt}
  &\hspace{-6pt}{(vy-u)/{(v-1)}}{\ < \hspace{.5pt}\frac{1}{8}\sqrt y}\,  -\frac{7}{8}
 \ ({\rm which~proves \  (P2(c))~using~(iv))}  
 \vspace{4pt}
 \\
\hspace{-6pt}&\hspace{-6pt}< &\hspace{-6pt} y/2 \hspace{.5pt} <\hspace{.5pt} y- 6 v 
\end{array}
\end{equation*}
 since $y>u/v >  12 v $  by (v) and (ii); and this proves  (P2(d)).
 \smallskip\smallskip

 (P4):  The relation $u=x+v(y-x)$ is just the present definition of $x$.
\qed

 \medskip   
{\em\noindent Remark.}
Condition (v) places $y$ in an interval of length roughly 
$ \frac{1}{8} \sqrt {{u/v}} $, which is fairly large by (ii).  
We still  need to verify the relatively obvious fact that this is large enough to make it possible to 
satisfy the remaining inequalities and congruences in the lemma.
 
 \begin{lemma}
   \label{Existence2}
   Given admissible integers  $v$ and $u$ such that {\rm$v\equiv 3$~(mod~12),}
   $v>16^5$
  and    $u\!>\! 800^2v^7 ,$  there are  integers $x,y_1,y_2$ and  $y$ behaving 
as stated in    {\rm (P1),}  {\rm (P2)}  and {\rm (P4).}
   \end{lemma} 
   
 \proof      
 We will use (i)-(v)  in the preceding lemma.   
 Let $u\equiv \delta$~(mod~4) with $\delta=\pm1$;  the remaining
 requirements in (i) and (ii) are among the present hypotheses.
  
  Since $v\equiv 3$~(mod~12) we can write  $v = 3 +12e+24m$  with $e\in\{0,1\}$.
  
  When $\delta = -1 $  let $y_2:=1$.  When $\delta = 1 $   let
 $y_2:=v+(6-e)(v-1)$,   so  $y_2 \equiv 1$~(mod~$v-1$)
 and  $y_2 \equiv 15$~(mod~$24$).   Clearly  $y_2<7v$.
 
We next define $y_1$.    Let $0<u'<v-1$ with  $u'  \equiv u$~(mod~$v-1$),    so $u'$ is odd. 
 
 Let $y_1:=u'+\Big(24Å\lceil u/\{24 v(v-1)y_2\}\big \rceil 
 +\frac{1}{2}\big((15-6e) +(23-6e) u' \big)\Big)(v-1)$.   
Since   $(15-6e) +(23-6e) u'$ is even
 we have   $y_1\equiv u'\equiv u  $~(mod~$v-1)$  and
 $y_1\equiv u'+ \big((15-6e) +(23-6e) u'\big)(1 +6e )  \equiv 15 $~(mod~$24)$.

 \emph{We claim that  $\,y_1,$ $y_2$ and  $y:=y_1y_2 $ behave as 
 required in} Lemma~\ref{Existence1}(iii-v).   
 
  \smallskip
   (iii):    $y=y_1y_2\equiv u\cdot1$~(mod~$v-1$).
  
  \smallskip
  (iv): Most of this is in our definitions of  $y_1$ and $y_2$,
  while  $u>24v(v-1)7v>24v(v-1)y_2$ 
  implies that $y_1>24(v-1)>7v>y_2$.
  
  \smallskip
  (v): 
For  the first part of (v),
   ${y=y_1y_2}> {\big (24 u  /\{24 v(v-1) y_2 \}\big)}(v-1) y_2 =u/v$.  
Next,  $u'<v-1, $  $\, y_2<7v$ and $ 800^2v^6< u/v$ imply that  \vspace{-2pt}
 $$
 \begin{array}{lll} 
 y=y_1y_2&\hspace{-6pt}< u'  y_2 +\big( u/v +24\cdot  1\cdot (v-1) y_2\big)
 +\frac{1}{2}(1  5+23 v)(v-1) y_2
 \vspace{3pt}
 \\
 &\hspace{-6pt}< (v-1)  7v +\big( u/v +24  (v-1)   7v \big)
 +\frac{1}{2}(1  5+23 v)(v-1)   7v 
 \vspace{2pt}
 \\
 &\hspace{-6pt} < u/v   +100 v^2(v-1) -1
 <u/v  +\frac{1}{8}\{(v-1)/v\}\sqrt {{u/v}}\, -1 
  \\
 &\hspace{-6pt}< u/v +u/v<u  ;
 \end{array} \vspace{-2pt}
 $$
  the ends of the last two lines take care of the remaining parts of (v).
 Now Lemma~\ref{Existence1} provides us with the required integer~$x$.\qed


 \subsubsection {Labeling}
 The structure of $  Y-X$  as  both a cyclic group and   a partial  STS 
   have nothing to do with one another, as observed by Moore \cite[p.~279]{Mo}. 
This independence is seen in  (M2) and (M3).
This allows us to label  the points of $Y-X$ in any way by the elements of $A$
using an arbitrary bijection $\pi\colon\!Y-X\to A$;
an element $y$ of $Y-X$ is labeled by $a:=y^\pi$,  which we abbreviate 
 by writing  $y=a$.

 \begin{lemma}
   \label{labeling}
The elements of $Y-X$ can be labeled by the elements of $A$ in such a way that$,$ 
if $k\in A_6,$  $\alpha\in {\rm Aut} A $  and 
the permutation $y\mapsto ky^\alpha$  of $A$
 is an automorphism of the partial Steiner triple system \ $Y-X,$  then $k=1$ and  $\alpha=1$. 
  \end{lemma}
 \proof 
By  (P1(a)) and  (P2(d)), $|Y-X| > 6\cdot 16^5$.  Then there are distinct points 
 $y_1,\dots,y_9$ of $Y-X$, and $x _0\in X$,  such that  
   the following are triples of $Y$:
 $$y_1,y_2,y_3 \qquad  y_3,y_4,y_5   \qquad x _0,y_6,y_7  \qquad  x _0,y_8,y_9.$$

 Let $c$ be a generator of $A$.  Label the $y_i$ using $A_6\cup \{c,-c,c^2\}$:
\begin{equation*}
  \begin{array}{lllllll}
 &y_1=1 &y_2=-1& y_3=c& 
  y_4=-c&y_5=c^2 ;   \rm \ \ and \    also
  \vspace{4pt}
  \\ 
&  y_6=\omega& y_7=-\omega&y_8=\omega^2&y_9=- \omega^2&
 \mbox{\rm if  some $ \omega \in A $  has  order  3.}&
   \end{array}
  \end{equation*}
  (The remaining points of $Y-X$ are 
  labeled by the remaining elements of $A$ in an arbitrary manner.
  Note that the points $y_6,y_7,y_8$ and $y_9$ are needed only if 
  $|A|$ is a multiple of 3.)
  Thus, we have the following triples in $Y$:
$$
1,-1,c \qquad  c,-c,c^2   \qquad x _0,\omega, -\omega  \qquad  x _0,\omega^2, -\omega^2 
$$
(where the last two are omitted if $|A|$ is not a multiple of 3).

 Now consider an automorphism $y\mapsto k y^\alpha$  of $Y-X$, where
$k\in A_6,$  $\alpha\in {\rm Aut} A $.  This sends the triple  $1,-1,c$ to $k,-k,kc^\alpha. $  
 If $k\in A_6-\{ \pm1\}$   then 
 this triple is $\omega^i,- \omega^i,k c^\alpha$ for 
  $i=1$ or 2;  but   the triple in $Y$ containing $\omega^i$  and  $- \omega^i$ is not contained in $Y-X$.
  
 Thus,  $k=\pm1$,  and $1,-1,c$ is sent  to $1,-1, k c^\alpha$,  so $k c^\alpha=c$.   
 
 If $k=-1$ then $c^\alpha=-c.   $  The triple $c,-c, c^2$ is sent to 
 $ -c^\alpha,  -(-c)^\alpha, -(c^2)^\alpha$,
  which is $c, -c, -(c^\alpha)^{2} $. Since  $-(c^\alpha)^{2}=-(-c)^{2}\ne c^{2}$, 
   this is impossible. 
   
 Then $k=1$, so $\alpha=1 $ since $c=c^\alpha$ generates $A$.
 \qed
     
  \subsection {Location of $PG(2,2)$ subsystems} 
  \label {Location of $PG(2,2)s$.} 
  We  need  structural properties of the STS $\,U$ defined in (P4).
  In this section we
   will not need any of the assumptions in    Section~\ref{Preliminaries}:
  only the definitions in (P4)  and the notation in (P5) are involved.
  
For $v\in V$ let  
$$ Y_v:=X\cup (v\times A) .$$
  By (M2) this is a subsystem of $U$
isomorphic to $Y$ (via the isomorphism $x\mapsto x,$ $ y\mapsto (v,y)$ for 
$x\in X,$ $ y\in A=Y-X$). 
   
  There are $PG(2,2)$ subsystems contained in $V\times 1$, 
  and  ones  contained in  $ Y_v \cong  Y$ for  $v\in V$.
  Another~possible type of  $PG(2,2) $ subsystem uses a triple $v_1,v_2,v_3$ in~$V$, $x\in X$,
  and elements $a_i\in A$:%
   \vspace{-1pt}
\begin{equation}
\label{oddball}
\begin{array}{lll}
\mbox{points: $x,$ $ (v_i,a_i)  ,$ $ (v_i,- a_i)  $  for $i=1,2,3 , $} 
\\\mbox{\hspace{34.5pt}for $x\in X,$ $a_i\in  A,$ $a_1a_2a_3=1$ and
 triples  $a_i ,-a_i,x$ in $ Y$} 

\vspace{2pt}\\
\mbox{triples: $ (v_i, a_i) ,  (v_i, -a_i) , x$ for $i=1,2,3$, and} 
\vspace{2pt}\\
\hspace{34.5pt}
\mbox{$ (v_1, \epsilon _1 a_1) , (v_2, \epsilon_2a_2) , (v_3, \epsilon_3a_3)$ whenever $\epsilon_i=\pm1$ and $\epsilon_1\epsilon_2\epsilon_3=1$.}
\end{array}
\end{equation}

\Remark
\label{oddball determined}
Two points $(v_1,a_1),   (v_2,a_2)$  with $v_1\ne v_2$ lie in 
at most one subsystem {\rm(\ref{oddball})}. \rm
 For, these points determine the  triple  $v_1,v_2,v_3$  and then all 
 $(v_i,\pm a_i)$.

 \defi\rm
 \label{V_{}}
Let  $S$ be a    subsystem of $ V$ and  $f \colon  \!S\to  A_6$ (cf.~(P5)) a function
such that $f({v_1})f({v_2})f({v_3})=1$ whenever $v_1,v_2,v_3$ is a triple in $S$.
Then 
$$V_{S,f}:=\{ (s,f(s))  \mid {s \in S}  \} $$
 is a subsystem of $U$,
and $V_{S,1}=S\times 1 \cong  V_{S,f}$ via $(s,1)\mapsto (s,f(s))$, using (M3): 
these subsystems are just variations on  the subsystem $V \times 1$.

      The  subsystems      $ Y_v$  and  $V_{S,f}$
 are   basic tools in our proof of Theorem~\ref{main},  and
  \begin{equation}
\label{int}
\mbox{$| Y_v\cap  V_{S, f}|\le 1$ for all $v,$ $S$ and   $f$.}
  \end{equation}
 \begin{lemma}
   \label{PG(2,2)s}  
   Every $PG(2,2)$ subsystem of $U$   either   is of type  {\rm(\ref{oddball}),} 
     lies in some $ Y_v,$ or has the form $ V_{S, f} $ 
  for a $PG(2,2)$  subsystem $S$  of $V$.
 \end{lemma}

 \proof
 If   a  $PG(2,2)$ subsystem $Z$ has the form  $\{(v_i,y_i)\mid 1\le i\le 7\}$ 
 with distinct $v_i$, then the 
 $v_i$ form an STS $ S$ by (M3);
  we may assume that the triples   in $Z$  are
 
$
 \vspace{1pt}
\hspace{21pt}
\begin{array}{llll}
(v_1,a_1) , (v_2,a_2) , (v_3,a_3)  && \mbox{so\, $a_1a_2a_3=1$}
\\
 (v_1,a_1) , (v_4,a_4)  ,(v_5,a_5)         & & \mbox{so\, $a_1a_4a_5=1$}   
 \hspace{130pt}
\\
(v_1,a_1) , (v_6,a_6)  ,(v_7,a_7)        && \mbox{so\,  $a_1a_6a_7=1$}
\\
 (v_3,a_3) , (v_5,a_5)  ,(v_7,a_7)       && \mbox{so\,  $a_3a_5a_7=1$}
\\
(v_3,a_3) , (v_4,a_4) , (v_6,a_6)  && \mbox{so\,  $a_3a_4a_6=1$}
\\
(v_2,a_2) , (v_4,a_4)  ,(v_7,a_7)  && \mbox{so\,    $a_2a_4a_7=1$}
\\
(v_2,a_2) , (v_5,a_5)  ,(v_6,a_6)  && \mbox{so\,  $a_2a_5a_6=1$}
\end{array}
$

 \vspace{1pt}
  \vspace{1pt}
\noindent
with $a_i\in A$.
 Multiplying  these equations, and also just the first three of them,  we find that  \vspace{2pt}
 $(\prod_ia_i)^3=1$ and $a_1^3 a_2a_3a_4 a_5a_6a_7=1$.   
 It follows that  $\prod_ia_i=\omega $ with $\omega^3=1$
  and $a_1^2=\omega^2$, so
 every $a_i ^6=1$.    Then $Z=V_{S,f}$ with $f(v_i):=a_i\in   A_6$.
 
If $Z$  contains a triple $(v,a_1), (v,a_2), (v,a_3)$
but does not lie in $ Y_v$, then it also contains a point  $(v_2,b_2)$ 
with $v_2\ne v$.  
By (M3),  if $v,v_2,v_3$ is a triple of $V$  then 
there are triples  
       $(v_2,b_2),(v,a_1),(v_3,b_3)$
and $(v_2,b_2),(v,a_2),(v_3,c_3)$
and $(v_2,b_2),(v,a_3),(v_3,d_3)$, 
and hence another triple  
$(v,a_1),(v_3,c_3),(v_3,d_3)$, which contradicts (M2).

Assume that $Z$  is not in any $ Y_v$.   
If $Z$ has  a triple $T\subseteq X$, then some point $(v,a)$ is in $Z$,
the triples joining $(v,a)$ to the points of $T$ all lie in both $Z$ and $Y_v$  by (M2),
and then $Z\subseteq Y_v$.
The only other possibility is that $Z$  is determined by three triples through 
some $x\in X$,  and hence contains triples
 
$\hspace{21pt}
\begin{array}{llll}
(v_1,a_1) , (v_2,a_2) , (v_3,a_3)  && \mbox{so $a_1a_2a_3=1$}
\\
x, (v_1,a_1) , (v_1,b_1)          && \mbox{so\,  $a_1, b_1,x$\,   is a triple in $ Y$}
\\
x, (v_2,a_2) , (v_2,b_2)           && \mbox{so\,  $a_2, b_2,x$\,   is a triple in $ Y$}
\\
x, (v_3,a_3) , (v_3,b_3)           && \mbox{so\,  $a_3, b_3,x$\,   is a triple in $ Y$}
\\
(v_1,a_1) , (v_2,b_2) , (v_3,b_3)  && \mbox{so\,  $a_1b_2b_3=1$}
\\
(v_1,b_1) , (v_2,b_2) , (v_3,a_3)  && \mbox{so\,  $b_1b_2a_3=1$}
\\
(v_1,b_1) , (v_2,a_2) , (v_3,b_3)  && \mbox{so\,  $b_1a_2b_3=1$.}
\end{array}
 \vspace{1pt}
  \vspace{1pt}
$

\noindent
The last three equations imply that $a_1a_2a_3b_1^2b_2^2b_3^2=1$,
so $1=(b_1b_2b_3)^2=(b_ia_i^{-1})^2$ and hence $b_i=\pm a_i$ for each $i$ (cf.~(P5)).
Since $ x, (v_i,a_i) , (v_i,b_i)   $ is a triple we  have $b_i=-a_i$
for each $i$,   so we are  in (\ref{oddball}).   \qed

\medskip 
For Moore   \cite[Sec.~10]{Mo},  the types of $PG(2,2)$ subsystems  of $\,U$
were isomorphism invariants of his STS
construction.\footnote{See Footnote 1.}
  He did not go into the detail involved 
in  (\ref{oddball})  or   a function ${f\colon\! S\to  A_6}$.%
 
\subsection{Finding critical subsystems} 
\label{Finding critical subsystems} 

The following key result  is based on  (P1(b)) and (P3(d)) 
 together with  ({\ref{int}).
\begin{prop}
\label{maximal} 

 Let $W$ be a  subsystem of $U$.
 \begin{enumerate}
 \item[(i)]
 If $W\cong  Y $ then  $W$  has the form    $ Y_v$  for some $v\in V$.

  \item[(ii)] 
  If $W\cong  V ,$  if $W$ meets every subsystem in 
 {\rm (i)}  in at most one point$,$ and if $W$ is disjoint from the intersection of the subsystems in  {\rm (i),}  then $W$ has the form $V_{V, f}$
  for some ${f\colon\! V\to  A_6}$.
   \end{enumerate}
  \end{prop}
   
\proof
(i)  Clearly $|W|=|Y|>|X|$.
 Let   $(v_1,a_1)\in W-X$.  
 We claim that $W \subseteq  Y_{v_1}$. 
 
 Assume that  $W\not\subseteq  Y_{v_1}$. 
 Let $(v_2,a_2)\in W$ with $v_1\ne v_2$.
Since $W\cong  Y $, by (P3(d))  
 there is a  $PG(3,2)$ subsystem   containing $(v_1,a_1) $ and 
 $ (v_2,a_2)$.  That subsystem  has three $PG(2,2)$ subsystems
 containing $(v_1,a_1)$ and $ (v_2,a_2)$, each of   
 type  $ V_{S, f}  $ or  as in 
 (\ref{oddball}), by Lemma~\ref{PG(2,2)s}; and  at least one has type 
  $ V_{S, f}  $ by  Remark~\ref{oddball determined}.    
In particular, $a_1\in A_6$.
Thus, $W\subseteq X\cup {(V  \times A_6)}$.~Then 
$| Y|=|W|\le |X|+6|V|$,   which contradicts  (P2(d)). \

\smallskip  
(ii) The stated intersection is $X$, so $W\subseteq  V\times  A$.  
  For   $v\in V$, $|W\cap Y_v|\le1$ implies that $v$ occurs in at most one pair $(v,a)  \in W $.   
    Since  $|W|=|V|$, it follows that  
 $W=\{(v,f(v)) \mid v\in V\} $ for some $f\colon \!V\to  A$.
 
 Since $W\cong  V $, by (P1(b))    every point of $W$ is in 
a $PG(2,2)$ subsystem, which by Lemma~\ref{PG(2,2)s} 
 has the form  $ V_{S, f'}  $   with   $|S|=7$ and  $f '\colon \! S\to  A_6$
 (since $|W\cap Y_v|\le1$ for each $v$).
Then $W\subseteq  V\times  A_6$, $f  $ maps to $  A_6$, and   by (M3) $f$
must behave as in Definition~\ref{V_{}}.
 \qed
 
 \subsection
{${\rm Aut}\hspace{.5pt}  U$  and   ${\rm Aut}  A$}
 \label{The occurrence of  $A$.}

Theorem~\ref{main} concerns   ${\rm Aut}  \hspace{.5pt}  U$:  
  
 \begin{prop}
  \label{Aut A}
${\rm Aut}\hspace{.5pt}  U\cong {\rm Aut}\hspace{.5pt}   V,$
and   ${\rm Aut} \hspace{.5pt} U$ leaves $V\times 1$ invariant$, $  inducing
${\rm Aut} ( V\times 1)\cong {\rm Aut} \hspace{.5pt} U$ on this subsystem
of $U$.
 \end{prop}
 \proof
 Let $h\in {\rm Aut}\hspace{.5pt}  U$.  We must show that $h$ is induced by some 
 element of  ${{\rm Aut} ( V\times 1)}\le {\rm Aut} \hspace{.5pt} U$
 (cf. (P4)).
 
  Proposition~\ref{maximal}(i) states that the subsystems $Y_v$ are uniquely determined for $U$. Then   Proposition~\ref{maximal}(ii) states that the subsystems $V_{V,f}$ are also 
  uniquely determined for $U$. 
It follows that $h$ sends $ V\times1=V_{V,1}$
to    $V_{ V, f}\subseteq V\times  A_6$
 for some $f\colon V\to A_6$,
and $h$ permutes the subsystems $ Y_v$.    

Since $h$ sends $X=\bigcap_{v\in V}Y_v$ to itself  it  also sends $U-X= V\times   A$ to itself.
In view of (M3), restricting $h$ to the first component in $ V\times   A$ 
induces an isomorphism $ \bar h\colon \!V\to   V$; by   (P4), 
\vspace{1pt}
$ \bar h$  is also induced by some 
$g\in {\rm Aut} ( V\times 1)\le {\rm Aut} \hspace{.5pt} U$.   Then
$ \bar h \bar g^{-1} =1$ on $V$.   We  will prove that $h=g$.
Replace $h$ by $hg^{-1}$,  so    $ \bar h=1$  on $V$.  
 {\em  The remainder of the proof consists of showing that $h=1$.}
 \vspace{2pt}

Since $  (v\times A)^h= v^{\bar h}\times A = v\times A $ and
 $( V\times 1)^h=V_{ V, f}$
we have  $(v,1)^h = (v,f(v))$  for all $v\in  V$.  

Since $h$ permutes the subsystems $ Y_v$, 
from $(v, 1),(v,1)^h= (v,f(v))\in   Y_v $ it follows that  $h$ leaves invariant
 every $ Y_v$.  
Let $(v,a)^h=(v,f_v(a))$ where $a,f_v(a)\in  A$.  
Then $(v,f_v(1))=(v,1)^h =  (v,f(v))   $.  Let  $b_v:= f_v(1) =f (v)\in A_6$.
 
{\em We will show that    $h$  acts   on $V\times  A$   by} 
 \begin{equation}
\label{va}
\mbox{$(v,a)^h=(v, b_va^\alpha)\,$ for all $v\in V, a\in A,$ and some $\alpha\in 
{\rm Aut} A.$}
 \end{equation}
 Let  $v_1,v_2,v_3$ be a triple   of $ V$.  Whenever $a_1a_2a_3=1$,  $a_i\in A$, 
 by (M3) we obtain
a  triple  $(v_1,a_1), $ $(v_2,a_2), $ $(v_3,a_3)$ and hence also 
its image under $h$:   the  triple 
 $(v_1,f_{v_1}(a_1)),$ $ (v_2,f_{v_2}(a_2)),$ $ (v_3,f_{v_3}(a_3))$, so 
 $f_{v_1}(a_1)f_{v_2}(a_2)f_{v_3}(a_3)=1$.  Then
 $$
 \mbox{$f_{v_1}(a_1)f_{v_2}(a_2)f_{v_3}(a_1^{-1}a_2^{-1})=1$   \  for all $a_1,a_2\in  A$.}
 $$
 Let $a_1=1$ and deduce that $f_{v_2}(a_2) =b_{v_1}^{-1}f_{v_3}(a_2^{-1})^{-1}$;
 while  $a_2=1$ yields $f_{v_1}(a_1) =b_{v_2}^{-1}f_{v_3}(a_1^{-1})^{-1}$.
 Then   $b_{v_1}b_{v_2}b_{v_3}=f_{v_1}( 1)f_{v_2}(1)f_{v_3}(1)=1$  and  (after replacing $a_i^{-1}$ by $a_i$)
 $$
b_{v_3} f_{v_3}(a_1  a_2  )= f_{v_3}(a_1     )  f_{v_3}(   a_2  ).
 $$
Now
$b_{v_3}^{-1} f_{v_3}(a_1  a_2  )= b_{v_3}^{-1}f_{v_3}(a_1     )  
b_{v_3}^{-1}f_{v_3}(   a_2  )$, so that
$f_{v_3}(a_1)=b_{v_3}a_1^\alpha$ for some $\alpha\in {\rm Aut}  A$
and all $a_1\in A$.
  Moreover, 
 ${f_{v_2}(a_2) =b_{v_1}^{-1}f_{v_3}(a_2^{-1})^{-1} = b  _{v_2}a_2^\alpha}$:
 we have the same automorphism $\alpha$ for all  $v\in V$.
 This proves (\ref{va}).
 
 By  (\ref{va})  and (M2),  if $v\in V$ then 
 $a\mapsto b_va^\alpha$ is an automorphism of 
 the partial Steiner triple system $Y-X$.  
 By (P6),  $\alpha=1$ and  $b_v=1$ for all $v$.  Then $h=1$ on $V\times A$.
 Since every point of $X$ is in a triple containing two points of $Y-X$, it follows that $h=1$, as claimed.  \qed
  
    \subsection
{Completion of proof}
\label{Completion of proof}
      
In (P1)-(P4) we provided the ingredients for the construction of an STS
$U$ using Section~\ref{XYV}.
Proposition~\ref{Aut A} determined Aut\hspace{.5pt}$U$.  

Moreover,  by  (P1(a))    and Proposition~\ref{Aut A},
$U$ has Aut\hspace{.5pt}$U$-invariant subsystems 
$V\times 1\supset V_{\scriptscriptstyle\star}\times 1$
such Aut\hspace{.5pt}$U \cong {\rm Aut}  
(V\times 1) \cong {\rm Aut}( V_{\scriptscriptstyle\star}\times 1)$
and Aut\hspace{.5pt}$U$ induces 
$ {\rm Aut}(V\times 1) $  and  $  
{\rm Aut}  ( V_{\scriptscriptstyle\star}\times 1)
 $ on the
respective subsystems.

Lemma~\ref{Existence2} states that we have dealt with all admissible  $u> 800 ^2v^7$. 
\qed
      
 \rema  \rm   
 \label{bounds}
{\bf Bounding  $N_{V_{\scriptscriptstyle\star}} $}. 
In Lemma~\ref{Existence1} we had $u> 800 ^2v^7$,  but
 the integer $v=|V|$ obtained in Corollary~\ref{=3}
   is much larger than $|V_{\scriptscriptstyle\star} |$.
  By  Proposition~\ref{VxD}(iv), 
  $|{V_{\scriptscriptstyle\star}}{\bm\times} D |$  is $O(|V_{\scriptscriptstyle\star} |^5)$, so 
$v$  is  $O(|V_{\scriptscriptstyle\star} |^{5\cdot 5})$. 
Thus,   $  N_{V_{\scriptscriptstyle\star}} $
is  $O(|V_{\scriptscriptstyle\star} |^{25\cdot 7})$.%
\medskip
  
  \noindent
  {\bf Bounding  $M_G $}. 
  In \cite{Me} and \cite{Ka}
   an STS\ $V_{\scriptscriptstyle\star}$ is constructed for which
  $G\cong{\rm Aut}V_{\scriptscriptstyle\star}$ 
 and $|V_{\scriptscriptstyle\star} |$ is $2^{O(|G|)}$.  By the 
  preceding paragraph  the same is true  for $M_G $.
  
    \em  
 This  bound for $M_G$  is ridiculously weak.  It seems
    likely that $M_G$  is polynomial in $|G|,$ but entirely new methods  
would be needed to prove that. \rm

    \rema  \rm 
    \label{comparison} 
 The argument  in \cite{Ca} depended on using pairs $X\subset Y$ 
provided by \cite{DW},  essentially for all possible  $x=|X|$ and $y=|Y|$
for which $y\ge  2x+1$.
The argument used here only had access to a more limited choice 
(P2(a)) of orders $y$  (cf. Lemma~\ref{direct product}(ii)).
In \cite{Ca} first $y-x$ was dealt with,  at which point
$x$ and $y$ were uniquely determined for given $v$ and $u$.  This approach can be used in our situation 
when $u\equiv -1$~(mod~4) but not when $u\equiv 1$~(mod~4).
Therefore we have started with a restricted choice of 
$y$, and  then $x$ is uniquely determined 
 for given $v$ and $u$ (Lemmas~\ref{Existence1}  and \ref{Existence2}).
Our problem was to have 
a suitably geometric $ Y$  of order $y$ with a subsystem   of the required 
order $x\le (y_0-1)/2< (y-1)/2$.
        
       \section {Partial Steiner triple systems} 
     \label {Proof PSTS} 
    \subsection {Theorem~\ref{PSTS}} 
      \label {Proof  of Theorem {PSTS}}
   We first note how our approach   differs from that of Cameron \cite{Ca}.  He observes:
``In the construction used to prove Theorem 1, if the subsystem
contains no triples, its automorphism group is the symmetric group $S_u$, while that of the embedding system is the general linear group $GL(u-1,2)$.''
In other words,  the PSTS  (partial Steiner triple system) might have too few triples.  
The first part of our proof  eliminates this possibility  (cf. Lemma~\ref{V'}(v)).

\definition  \rm
\label{paths}
 Let the  PSTS  $Q_k(x), k\ge 2,$ have the following triples (using two ``paths'' of $k$ triples in the first two rows and an additional point $2k+1$):
 $$
 \begin{array}{lllllllllll}
  x,1,2   &2,3,4\hspace{28.45pt} \  4,5,6\hspace{18.9pt}\dots \quad2k-2,2k-1,2k
\\
x',1',2'   &2',3',4' \qquad  \   4',5',6'\quad\dots \quad(2k-2)',(2k-1)',(2k)'
\\
 2k+1,x,x'\ &  2k+1,i,i'  \mbox{\hspace{2.4pt}  \    for $1\le i\le 2k$}  \ \ 
\\
x,3 ,(2k)'  & x,3',2k  \end{array}
 $$
 
\Remark\rm
\label{paths2}
The following properties of $Q_k(x)$ are straightforward:
\begin{enumerate}
\item[(1)] $Q_k(x)$ has $4k+3$ points, 
\item[(2)] every point   is in at least two  triples$ ,$
\item[(3)]  the point
$2k+1$ is in $2k+1\ge 5$ triples, 
$x$ is in four triples and every other point is in at most  three triples,
\item[(4)]every point is in the  union of the triples containing $2k+1$,   
and
\item[(5)]${\rm Aut}Q_k(x)=1$.
\end{enumerate}
(For (5),  an automorphism must fix $x$ and $2k+1$, then also 
$x'$, $1'$, $1$,   $2$, \dots.)
\smallskip

Let $V$ be an $n$-point PSTS as in Theorem~\ref{PSTS}.  We may assume that $n\ge2$.

\begin{lemma}
\label{V'}
There is a PSTS\, $V'$ such that
\begin{enumerate}
\item[(i)]   
$ V $ is an  ${\rm Aut} \hspace{.5pt}  V'$-invariant   subsystem  of $V',$
\item[(ii)]${\rm Aut}\hspace{.5pt} V'$ induces 
${\rm Aut}\hspace{.5pt}  V$ on $V,  $
\item[(iii)]${\rm Aut} \hspace{.5pt} V'\cong{\rm Aut}\hspace{.5pt}  V,$   
\item[(iv)] $n':=|V'|\ge 22,$  and
\item[(v)]
every point of $V'$ is in at least two triples of $V'.$
\end{enumerate} 
\end{lemma}
\proof
 For every point $x$ of $V$,  attach $Q_{n }(x)$ to $V $  so that $V \cap Q_{n}(x) = x$  and the  
 $n$ PSTSs $Q_{n}(x)$ are pairwise disjoint. The union of $V$ and these 
 PSTSs (also using the union of  their sets of triples) is a new PSTS   $\,V' $ having $n'$  points, where 
 $n' = n|Q_{n}(x)| = n(4n +3) \ge 22 $,  which proves (iv).
 
 Condition  (i) is clear, (v) holds in $V'$ by Remark~\ref{paths2}(2),
 and  (ii) follows from the fact that  all $Q_{n }(x)$ are isomorphic  and are pairwise disjoint.  
  
  It remains to prove (iii).  By Remark~\ref{paths2}(3),
 any subsystem $Q$ of $V'$ isomorphic to $Q_n(x)$ has a point $z$ in $2n+1$ triples  of $Q$.  Since
\vspace{2pt}
 $V'=\bigcup_{x\in V}Q_n(x)$,   again by Remark~\ref{paths2}(3)
each point of $V'$ is either in $2n+1$ triples of $V'$,
at most 3 triples, or (for points of $V$) between 4 and $4+(n-1)/2 <2n+1$ triples.
Then $z\notin V$ and $z\in Q_n(x)$ for a unique $x$.
By Remark~\ref{paths2}(4),  the union of the triples of $V'$ containing $z$ is both $Q$ and $Q_n(x)$, so  $Q=Q_n(x)$. 

 Thus, $V'$ determines the points of $V$,  
  any element of ${\rm Aut }  V'$ induces an element of ${\rm Aut }  V,$ and this yields a homomorphism from ${\rm Aut }  V'$ onto ${\rm Aut }  V.$ Its kernel fixes every point $x$ of $V$, and hence is 1
by Remark~\ref{paths2}(5). \qed
 
      \smallskip \smallskip
In the rest of the proof we ignore $V$ and work only with  $V'$.
By  Theorem~\ref{main},  \emph{it suffices to construct one STS\, $U$ 
having $ V '$ as an  ${\rm Aut} \hspace{.5pt}U$-invariant   subsystem  
such that ${\rm Aut} \hspace{.5pt}U\cong{\rm Aut}\hspace{.5pt}  {V' }$
and  ${\rm Aut}\hspace{.5pt} U$  induces ${\rm Aut}\hspace{.5pt}  V'$  on $V'$.}

 As in \cite{Ca}, we use the projective space $P=PG(n'-1,2)$
 whose points are the $2^{n'} - 1$ nonempty subsets of (the set of points of)   $V' $,
 the lines of $P$ being all triples of subsets   of $ V'$ whose symmetric difference is empty. Any permutation of the points of $V'$ extends uniquely to an automorphism of $P.$   Every  point $w$ of  $P$ has   size $|w|$ as a subset of~$V'$.

 Again as in \cite{Ca},  we construct from the STS $P$ and the PSTS $V'$
  an STS $U$  whose points are those of $P$, as follows: for every triple 
 $a, b, c $ of  $V'$,   replace the triples
\begin{equation}
ab, ac, bc \qquad a,b,ab \qquad a,c,ac     
\qquad b,c,bc
\end{equation}
of $P$ by the  new triples
\begin{equation}
\label{triples in S}
a, b, c    \qquad a,ab,ac\qquad b,ab,bc   
\qquad c,ac,bc
\end{equation}
(by abuse of notation, we write $a$ and $ab$ for $\{a\}$ and $\{a, b\}$, respectively). This produces a new STS $\,U$, because the new triples cover exactly the same pairs of points as the old ones. 

Note that
\begin{equation}
\label{2}
  \mbox{\em Every point $ab$ is in at most two triples of $U$ that are not lines of $P,$} 
\end{equation}and that ${\rm Aut }  V'$ induces a subgroup of ${\rm Aut }  \hspace{.5pt}U$ (as in \cite[p.~468]{Ca}).
 Moreover, 
 \begin{equation}
\label{block is line}   
 \mbox{\em A line of $P$ is also a triple of $U$ if it contains a point $w$ with $|w|>2$.} 
 \end{equation}
 
\begin{lemma}
\label{triples}
The lines of $P$ can be determined using the triples in $U$.
\end{lemma}
 \proof
We will recover  the line
 $\langle x,y\rangle$ of $P$  determined by any given   distinct points $x,y$ of $U$.   
For every point $p $ of $U $   not in the triple of $U$ determined by $x$  and  $y$
there are distinct triples 
$$p,x,x_1 \qquad p,y,y_1 \qquad x_1,y_1,z    \qquad p,z,q$$
of $U$, producing a 7-set  $U(p,x,y) :=\{p,x,y,x_1,y_1,z,q\}$ of points of  $U$.   

There are at least $(2^{n'-2}-1) -n'-{\binom{n'}{2}}$ planes 
\vspace{3pt}
of $P$ containing $x$ and $y$ but containing no  point $w$ of $P$ with $|w|\le 2$.
Every point   $w\notin \langle x,y\rangle$ in such a plane has $|w|>2$;
by   (\ref{block is line}),  every such  plane has the form $U(p,x,y)$ for any of its four points
$p\notin \langle x,y\rangle$. Thus, by Lemma~\ref{V'}(iv),  if $p$ is one of at least 
$4(2^{n'-2}-1-n'-{\binom{n'}{2}})>\frac{3}{4}|U|$ 
points in the union $S $ of these planes but not in    $\langle x,y\rangle$, then

\begin{itemize}
 \item [\rm(i)]    every set $U(p,x,y)$   occurs for at most four points $p$,
 and
  \item[\rm(ii)] distinct sets $U(p,x,y)$  have the same intersection  of size 3.
\end{itemize}
(The intersection in (ii) is the line  $\langle x,y\rangle=\{x,y,z\}$.)

If $S '$  is another set of more that  $\frac{3}{4}|U|$  points  satisfying (i)-(ii), then 
$|S\cap S' | \ge\frac{3}{4}|U|+\frac{3}{4}|U|- |U|=\frac{1}{2}|U|=\frac{1}{2}(2^{n'}-1)  
\ge \frac{1}{2}(2^{22}-1) $,  and hence by (i)  $
S\cap S'$ contains distinct sets $U(p,x,y)$.
Those sets produce  the same set  of size 3 in  (ii).
Thus, we have obtained the line  $ \langle x,y\rangle$  of $P$ using the triples of $U$.   \qed  
\smallskip\smallskip

Note that we have not yet used Lemma~\ref{V'}(v).
 
\begin{prop}
${\rm Aut}\hspace{.5pt} U\cong {\rm Aut} V'\cong {\rm Aut} V.$
\end{prop}

\proof  
We now have the triples in $U$ and the triples in $P$.
 Let $T$ denote the set of triples of $U$ that are not triples in   $P$
(these are the triples in (\ref{triples in S}), 
and hence consist of points such as $a\in V'$ or $ab$).  Let $x \in U$.
\begin{itemize}
 \item [\rm(1)]  If  $x$ is  in at least four triples in $T,$ then $x\in V'$ by  (\ref{2}).
 \item[\rm(2)]  Every point of $V'$ is in at least four triples in
 $T,$ by  Lemma~\ref{V'}(v)  and  (\ref{triples in S}).
 
\end{itemize}
  Thus {\em$U$ uniquely determines $V',$}   so ${\rm Aut } \hspace{.5pt} U$ induces a subgroup of 
  ${\rm Aut }\hspace{.5pt}   V'.$ Since ${\rm Aut }  V'$ induces a subgroup of 
  ${\rm Aut }  \hspace{.5pt}U,$
  by  Lemma~\ref{V'}(iii)       we have
  ${\rm Aut }\hspace{.5pt}  U \cong {\rm Aut } \hspace{.5pt}  V' \cong  {\rm Aut }  V,$  and we are done. \qed
  \medskip

{\noindent\em Proof of {\rm  Theorem~\ref{PSTS}}.}
We have embedded the original PSTS $V$ into an STS $\,U$ such that
${\rm Aut }\hspace{.5pt}  U $ leaves $V$ invariant, induces
${\rm Aut }  V $  on $V$, and is isomorphic to ${\rm Aut }  V $.
Now apply  Theorem~\ref{main} to $U$ (in place of $V$) in order to obtain  STSs 
behaving as in Theorem~\ref{PSTS}}.   \qed
    
\subsection{Corollaries} \label{corollaries}
 We note some   consequences of Theorem~\ref{PSTS}.
 We will use   a natural  graph on  the points of  a   PSTS $W$, 
with two points joined if they are in a triple.  If this graph is not connected we can embed $W$ in an arbitrarily large STS, which is clearly connected (preservation of the automorphism group is even possible by Theorem~\ref{PSTS}, but this will not be needed). 

\begin{coro}
\label{2 STS}
Given partial Steiner triple systems $V$ and $W,$  there is an integer $N_{V,W}$ such that$,$
for  each admissible   $u\ge N_{V,W} ,$   there is a Steiner triple  system $U$ on $u$ points having 
a subsystem $W'\cong W$ and 
  an  ${\rm Aut}\hspace{.5pt} U$-invariant   subsystem    $V'\cong V$   
  with $ W' \cap V'=\emptyset$  such that  ${\rm Aut}\hspace{.5pt} U   \cong{\rm Aut} V'$  and
   ${\rm Aut}\hspace{.5pt} U  $  induces  ${\rm Aut} V'$  on~$V'$.
\end{coro}

\proof
Here $V$ and the desired $U$ are as usual, the new aspect is to include $W$ as well;
we have no  control over the PSTS $V\dot \cup W$.
By the preceding remarks, we may assume that   $W$ is an STS and hence is
 connected, and that $W\cap V=\emptyset$,
 $n=|W|\ge2$   and $n>|V|$. Let $x_1,\dots,x_n$ be the points of $W$.
Then Definition~\ref{paths} applies with $k=i+n\ge 3$;
attach pairwise disjoint PSTSs  $Q_{i+n}(x_i) $ to $W$ in such a way that   
  $Q_{i+n}(x_i) \cap W =x_i$ for every $i$.
 The union of $W$ and all $Q_{i+n} (x_i) $  is  a connected PSTS $\widehat W$.
 
 Every $Q_{i+n}(x_i)$ has a unique point in  $2(i+n)+1\ge 2n+3 >\frac{1}{2}(n-1)+4$ triples 
 (by Remark~\ref{paths2}(3)),
 and  $\widehat W$ has no other such points ($x_i$ is in at most 
 $\frac{1}{2}(n-1)+4$ triples of $\widehat W$,   again by Remark~\ref{paths2}(3)).
Then $W$ can be recovered from  $\widehat W$  using   Remark~\ref{paths2}(3)-(4). The PSTSs 
\vspace{1pt} 
$Q_{i+n}(x_i), 1\le i\le n, $ have different orders,    so
from Remark~\ref{paths2}(5)   it follows that  ${\rm Aut} \widehat W=1$. 
 
Since $|\widehat W|>|V|$  and  $\widehat W$  is a connected component  of the graph on the disjoint union $\widehat W\dot{\cup} V'$ of $\widehat W $  
and  $ V'\cong V$,\, 
 ${\rm Aut} (\widehat W\dot{\cup} V')$ leaves $\widehat W$
invariant and hence acts on $V'$.  Then ${\rm Aut} (\widehat W\dot{\cup} V')\cong  {\rm Aut} V'$,  so
  apply Theorem~\ref{PSTS} to  $\,\widehat W\dot{\cup} V'$.\qed

\medskip
The first step in the above proof was to embed an arbitrary STS into one 
whose   automorphism group is trivial.  This suggests a strengthening of Theorem~\ref{corollary}:
 
\begin{coro}
\label{corollary2}
 If   $V _1,\dots,V_m $ are partial Steiner triple systems and $G$ is a finite group$,$
 then there is an integer $N_{V _1,\dots,V_m,G}$  such that$,$ 
for  each  admissible $u\ge N_{V _1,\dots,V_m,G},$ 
   there is a Steiner triple system $U$ on $u$ points
such that $V _1,\dots,V_m $ are   isomorphic to pairwise disjoint
subsystems of $U$  and ${\rm Aut} \hspace{.5pt}U\cong G.$

 \end{coro} 
\proof
Let $V$ be an STS with    ${\rm Aut} \hspace{.5pt}V\cong G $
\cite{Me}.
Apply the preceding corollary to 
 $V$ and    $W$, where $W$ is the disjoint union of (copies of)  $V _1,\dots,V_m $. \qed
 
\begin{coro}
\label{2 groups}
If $G$ and $H$ are finite groups  then there is an integer 
$N _{G,H}$  such that$,$ 
for  each  admissible   $u\ge N _{G,H},$   there is a Steiner triple system $U$ on $u$ points having a subsystem $W$ such that ${\rm Aut}\hspace{.5pt} U\cong G $ and ${\rm Aut} \hspace{.5pt} W\cong H $.
\end{coro} 
\proof
Let $V$ and $W$ be  STSs such that ${\rm Aut} \hspace{.5pt} V\cong G $  and
${\rm Aut} W\cong H $ \cite{Me}.   Apply  Corollary~\ref{2 STS}
 to the pair $V,W$ in order to obtain an STS $\,U$ behaving as stated.  \qed

\medskip
This corollary can be iterated in two ways:  one involves disjoint 
subsystems with arbitrary given automorphism groups; another involves
a nested sequence of subsystems with 
arbitrary given automorphism groups.

Our final corollary concerns retaining a {\em subgroup} of the automorphism group of an STS 
but not the full automorphism group.
Notation: If $G$ is a group acting on a set $X$, and if $Y\subset X$, then the
\emph{set-stabilizer}   $G_Y$ is $\{g\in G  \mid \mbox{$g$ sends $Y $ to itself\}}$.
(In the corollary $V_1$ need not be ${\rm Aut} \hspace{.5pt} V $-invariant, and
${\rm Aut}\hspace{.5pt}  V_1 $ need not be a subgroup of  ${\rm Aut} \hspace{.5pt} V $.)

\begin{coro}
\label{subSTS}
If   $V_1 $ is a subsystem of  order $>1$ of
a partial Steiner triple system $V ,$ then there is an integer  $N'_{V,V_1}\!$  such that$,$ 
for  each  admissible   $u\ge N'_{V,V_1},$ there is a Steiner triple system $U$ on $u$ points
having $V $  and $V_1$    as   ${\rm Aut} \hspace{.5pt}U$-invariant   subsystems   
such that  ${\rm Aut}\hspace{.5pt} U \cong  ({\rm Aut}\hspace{.5pt} V)_{V_1}$ 
 and ${\rm Aut}\hspace{.5pt} U$   acts on $V$ as  
  $({\rm Aut}\hspace{.5pt} V)_{V_1}$.
 
 \end{coro} 
\proof
First we replace $V$ by an STS:  use Theorem~\ref{PSTS} to find an  STS  $\,\hat V $ containing 
$V$  such that  $  {\rm Aut}\hat V$ leaves $V$ invariant, induces  ${\rm Aut} V$ on $V$ and is isomorphic to  ${\rm Aut} V$. Then
$ ({\rm Aut}\hat V)_{V_1} \cong  ({\rm Aut}\hspace{.5pt} V)_{V_1}$.

Let $z$ be a new point, and let $x\mapsto x'$ be a bijection from $V_1$ to a set $V_1'$
disjoint from $\hat V\cup \{ z\}; $ this bijection turns $V_1'$ into  a PSTS.
Form a   PSTS\,  $W, $ with $\hat V\cup \{ z\}\cup V_1'  $ as its set of  points, by
using the triples in $\hat V\cup V_1'  $  and  
\vspace{1pt}
including  a new triple  $x,z,x'$ for every $x\in V_1$. 
Every $g$  in   $ ({\rm Aut}\hat V)_{V_1}$ acts as an automorphism of $W$ via  $z^g=z$ and 
$(x')^g = (x^g)'$  for $x\in V_1$.
 
 The set $V_1$ is uniquely determined as the set  of points  of $W$
 lying in triples~with the maximal number  of other points (namely,  $(|\hat V|-1)+2$ points); 
 $\hat V -V_1$ is uniquely determined as the set of points of $W$ lying in  triples  with   exactly  $|\hat V| -1$ points.
 Then ${\rm Aut}\hspace{.5pt} W$ induces 
 $ ({\rm Aut} \hspace{.5pt} \hat V)_{V_1}$ on $\hat V$.
 
 Let $K$ denote  the  pointwise stabilizer of $\hat V$ in 
 ${\rm Aut} \hspace{.5pt}  W$.  For distinct $x,y\in V_1$, 
 the triples $x,z,x'$ and $y,z,y'$ meet at $z$,  so $K$ fixes $z$ and then also all points of $V_1'$.
Thus, $K=1$ and  ${\rm Aut}\hspace{.5pt} W=({\rm Aut}\hspace{.5pt}  \hat V)_{V_1}$.
 Now apply Theorem~\ref{PSTS} to $W$. \qed

\Remark
If $G$ and $H$ are finite groups with $G>H,$  then there is an integer 
$N '_{G,H}$  such that$,$ for  each  admissible   $u\ge N _{G,H},$ 
there  a Steiner triple system $U$   on~$u$ points having a subsystem $W$
 such that ${\rm Aut} \hspace{.5pt}U\cong G $ and
 $({\rm Aut}\hspace{.5pt} U)_W \cong {\rm Aut} W\cong H.$
\rm 
The proof involves a few straightforward changes in \cite[Sec.~2]{Ka}, which we  briefly outline.
     
{\bf  1.}  Let  $\Gamma$  be  an $n$-vertex connected graph   having a connected
   induced subgraph  $\Gamma'$ such that
  ${\rm Aut} \Gamma  \cong G$,  \,$G$ acts semiregularly on the vertices
  of $\Gamma$, and
$({\rm Aut} \Gamma)_{\Gamma'} \cong {\rm Aut} \Gamma'  \cong H$.   (Use the standard colored Cayley graphs for $G$ and $H$  and replace colored edges by  suitable   graphs.)
  
 {\bf  2.}    
 Consider the vector space $V$ in Section~\ref{Increasing $V$.}.  
 In order to conform to the notation in  \cite[Sec.~2]{Ka} let  $V_F$ and $V$ denote this as an 
 $F$-space and as an  ${\bf F}_2$-space, respectively.
  Assume that $G$ acts on the basis  of $V_F$  as it does  on the vertices of   $\Gamma$. 
  Let $V'_F$ denote  the $F$-span of the vertices of $\Gamma'$, and let $V'$ be 
  $V'_F$ viewed as an ${\bf F}_2$-space.

 {\bf  3.}     In \cite[Sec.~2]{Ka} there is a construction of an STS $\,U$  
 with   ${\rm Aut} \hspace{.5pt} U \cong G$,
  using $V_F$, $V$  and   $\Gamma$, and two auxiliary STSs on 15 points.    
  Restricting the construction to $V_F'$, $V'$  and   $\Gamma'$
    produces a subsystem $U'$ of $U$  obtained
    using these ingredients  in the same manner that  $U$  was.  In particular, 
      $({\rm Aut} \hspace{.5pt} U)_{U'} \cong {\rm Aut}  \hspace{.5pt}U'  \cong  H$, as required. 
  
            \smallskip   \smallskip   \smallskip
{\noindent\em Acknowledgements.} \rm
We are grateful to a   perspicacious, precise and patient referee,
whose comments led to significant improvements.
The research of the second author was supported in part by a grant from the Simons Foundation.


\end{document}